\documentclass[12pt]{amsproc}
\usepackage{amsmath,amssymb,mathrsfs}
\newcommand\ad{\operatorname{ad}}

\title[Darboux and a noncommutative bispectral problem]{The Darboux process and a noncommutative bispectral problem} 
\author{F. Alberto Gr\"unbaum}
\thanks{The author was supported in part by NSF Grant \# DMS 0603901}
\date{}
\address{Department of Mathematics \\ University of California \\ Berkeley,
CA\ \ 94720}

\begin{document}

\begin{abstract}
The Darboux process, also known by many other names, played
a very important role in some extremely enjoyable joint work that Hans and I did 25 years ago. I revisit a version of this problem in a case when scalars are replaced by matrices, i.e., elements of a non-commutative ring.
Many of the issues studied here can be pushed to the case of a ring with identity, but my emphasis is on very concrete
examples involving $2 \times 2$ matrices.
\end{abstract}

\maketitle

\begin{center}

The material in this paper was presented at a conference in honor of Hans Duistermaat, Aug 20-24 2007.
\end{center}
\section{The bispectral problem}

About 30 years ago I was extremely lucky.
I gave a talk in Berkeley where I mentioned
the following problem:

Find all nontrivial instances where a function $\varphi(x,k)$ satisfies
\[
L\left( x,\frac {d}{dx} \right)\varphi(x,k) \equiv (-D^2 + V(x))\varphi(x,k)
= k\varphi(x,k)
\]
as well as
\[
B\left(k,\frac {d}{dk} \right)\varphi(x,k) \equiv \left( \sum_{i=0}^M
b_i(k)\left( \frac {d}{dk} \right)^i\right) \varphi(x,k) =
\Theta(x)\varphi(x,k).
\]
All the functions $V(x),b_i(k),\Theta(x)$ are, in principle, 
arbitrary except for
smoothness assumptions. Notice that here $M$ is arbitrary (finite). 

I was fortunate that Hans was in the audience, and about a week later he came up with a tool to attack this problem. After a few months of intense work, mainly by
slow mail, we found ourselves with a rather nice picture.

The complete answer to the problem
is given as follows:

\bigskip
\noindent
{\bf Theorem 1 \cite{DG}.} {\em If $M = 2$, then $V(x)$ is (except for
translation) either $c/x^2$ or $ax$, i.e., we have a Bessel or an Airy case.
If $M > 2$, there are two families of solutions.}

\begin{itemize}
\item[a)] {\em $L$ is obtained from $L_0 = -D^2$ by a finite number of
Darboux transformations $(L = AA^* \rightarrow {\tilde L} = A^*A)$.  In this
case $V$ is a rational solution of the Korteweg-deVries equation and all rational solutions of KdV decaying at
infinity show up in this fashion.}

\item[b)] {\em $L$ is obtained from $L_0 = -D^2 + \frac {1}{4x^2}$ after a
finite number of rational Darboux transformations.}
\end{itemize}

\bigskip
It was later observed in \cite{MZ} that in the second case we are dealing
with rational solutions of the Virasoro or master symmetries of KdV.

In the first case the space of common solutions has dimension one, in the
second it has dimension two. One refers to these as the rank one and rank two
situations. In \cite{DG} one finds several other equivalent descriptions of the
solution such as those in terms of the monodromy group of the equation.

Observe that the ``trivial cases'' when $M=2$ are self-dual in the sense that
since the eigenfunctions $f(x,k)$ are functions either of the product $xk$ or
of the sum $x+k$, one gets $B$ by replacing $k$ for $x$ in $L$.

My reasons for asking this question can be traced back to an effort 
to understand some work on ``time-and-band-limiting'' that had led me to isolate certain properties of well known special functions. For an example relating to
orthogonal polynomials see \cite{G1}. For a more updated versions of this 
connection between the bispectral problem and the issue of time-and-band-limiting see, \cite{G2,G3,GY2}.

The work with Hans gave rise to a large number of papers by other people,
some of which can be found in the list at the end of this paper. This listing is far from complete and I apologize for the omissions.

It may be appropriate to observe that what we are calling the Darboux process has been
reinvented many times, including in the work of some rather well known people, see for instance \cite{Sc,IH}. Reference \cite{YZ} talks about the Geronimus
transformation, from 1940, and its inverse the Christoffel transform. It is clear that the first one (as noticed in \cite{YZ}) has a lot in common with what
we are calling the Darboux transformation. See also \cite{SVZ,Z}.

\section{The Bochner Krall problem}

In a series of papers with Luc Haine, \cite{GH1,GH2,GH3} and then with Luc and Emil Horozov, see \cite{GHH1,GHH2},
we noticed that a large class of polynomials  $$p_n(k)$$ that satisfy three term recursion relations
in the variable $n$, as well as differential equations in the variable $k$ can be obtained 
by an application of a similar Darboux transformation starting from the so called classical
orthogonal polynomials of Jacobi, Laguerre and Hermite.
In this case one goes from a tridiagonal matrix $$L_0=AB$$ to a new tridiagonal
matrix $L$ by factorizing the first one (or a function of it) as a product of two bidiagonal
matrices. As indicated in \cite{GH1} some form of this method is given in \cite{MS1}.

In this case one runs into the
Toda flows and its master symmetries. Further work on these lines can be found in \cite{GY1} and \cite{GH4,GH5} and for a very nice survey of all this material see \cite {H2}.

The origins of this line of work is contained in papers such as \cite{Bo,Koo,Kra1}.

\section{A matrix valued version of the Darboux process for a difference operator}

Consider the block tridiagonal matrix $L_0$
\[
L_0 = \begin{pmatrix}
B_0 & I  \\
A_1 & B_1 & I  \\
&\ddots &\ddots &\ddots
\end{pmatrix}
\]
where all the matrices $A_i,B_i$ are of size $N \times N$ and $I$ denotes the $N
\times N$ identity matrix.

If we try to factorize this in the form  $$L_0= \alpha  \beta$$ where
\[
\alpha = \begin{pmatrix}
\alpha_0 & I  \\
0 & \alpha_1 &  I  \\
&\ddots &\ddots &\ddots
\end{pmatrix}
\]
and 
\[
\beta = \begin{pmatrix}
I  \\
\beta_1 &  I  \\
&\ddots &\ddots 
\end{pmatrix}
\]
with all the matrices $\alpha_i,\beta_i$ are of size $N \times N$, and then 
define the matrix $$L= \beta \alpha$$  we have that
\[
L = \begin{pmatrix}
{\tilde B_0} & I  \\
{\tilde A_1} & {\tilde B_1} & I  \\
&\ddots &\ddots &\ddots
\end{pmatrix}
\]
where all the matrices ${\tilde A_i},{\tilde B_i}$ are of size $N \times N$, and the following relations
hold
$$\beta_n = B_{n-1} -\alpha_{n-1} , \alpha_n= A_n \beta_{n}^{-1} $$  
which then gives
$${\tilde B_n}= B_n-\beta_{n+1}+\beta_n  = B_{n-1}- \alpha_{n-1}+\alpha_n$$
and  
$$ {\tilde A_n}= \beta_n A_{n-1} \beta_{n-1}^{-1}  = \alpha_n^{-1} A_n \alpha_{n-1}.$$ 
These expressions are valid for $n=2,3,\dots$ in the case of ${\tilde A_n}$ and for $n=1,2,\dots$ in the
case of ${\tilde B_n}$.

Above we take $\beta_0=0$, so that ${\tilde B_0}=\alpha_0=B_0-\beta_1$.  We also need to take
${\tilde A_1}=(B_0-\alpha_0) \alpha_0$.

A moment's thought gives that once $L_0$ is given, the only free parameter is the 
{\em matrix} $\alpha_0$.

Just as in \cite{GH1}, and in spite of the fact that one is dealing here with a 
semi-infinite block tridiagonal matrix, it is possible to see the connection between this
construction and that in \cite{MS1}. One puts $$\beta_n= -\phi_n \phi_{n-1}^{-1}$$ and then 
notices that this amounts to choosing $\phi$ in the null-space of $L$.  Since $L$ is not
doubly-infinite we seem to have lost some freedom in picking this subspace, but this can
be remedied  as in \cite{GH1} by considering $L$ as a limit of an appropriate
doubly-infinite matrix with a rich null-space.

\section{Fancier versions of the Darboux process}

It is well known that
it is useful to extend the original method of Darboux
consisting of going from $$L_0=AB$$ to $$L=BA$$ in an appropriate way.
 Notice that in the standard case we have
$$B L_0=L B.$$

The new idea is to allow for an arbitrary banded matrix ( or a differential operator)  $U$ and to declare $L$ a Darboux
transform of $L_0$ as long as we have $$U L_0=L U.$$

In several of the uses of the original Darboux's method one needs to apply it repeteadly
and this fancier version of the method takes care of that.

One should also keep in mind the results in \cite{GHH1,GHH2,H2} where the usual factorization followed by a reversal of the factors is applied not directly to $L$ but to a 
constant coefficient polynomial in $L$.

I thank Jose Liberati for pointing out to me that at the very end in \cite{GGRW} one finds an application
of the theory of quasideterminants (a notion that goes back to Cayley) to obtain expressions for matrix valued orthogonal polynomials in terms of their matrix valued
moments. Many of these results, as well as others, have been derived independently by making good use of the notion of Schur
complements in L. Miranian's Berkeley thesis. The main new results are contained in \cite{M}.
In the next section we give a very brief look into the theory of matrix valued orthogonal
polynomials and a short guide to the literature that is relevant to us.

One should remark that this same theory of quasideterminants has been studied in 
connection with a certain Darboux process for a matrix Schroedinger equation in \cite{GV}. In this case
one needs to consider this fancier version. For a nice use of quasideterminants in our context see \cite{BL}.

The matrix version of the Darboux process for a difference operator discussed in the
previous section could be extended in this fancier fashion too.

\section{Matrix valued orthogonal polynomials}

Given a self-adjoint
positive definite matrix valued weight function $W(t)$, Krein, see \cite{K1,K2}, 
considers the skew symmetric bilinear form defined for any pair of
matrix valued functions $P(t)$ and $Q(t)$ by the numerical matrix
$$
\langle P,Q\rangle=\langle
P,Q\rangle_W=\int_{\mathbb{R}}P(t)W(t)Q^*(t)dt,
$$
where $Q^*(t)$ denotes the conjugate transpose of $Q(t)$.

Proceeding as in the case of a scalar valued inner product Krein
proves that there exists a sequence $(P_n)_n$ of matrix
polynomials, orthogonal with respect to $W$, with $P_n$ of degree
$n$ and monic.

Krein goes on to prove that any sequence of monic orthogonal
matrix valued polynomials $(P_n)_n$ satisfies a three term
recurrence relation
\begin{equation}
A_nP_{n-1}(t)+B_nP_n(t)+P_{n+1}(t)=tP_n(t),
\end{equation}
where $P_{-1}$ is the zero matrix and $P_{0}$ is the identity matrix.
These coefficient matrices enjoy certain properties: in particular
the $A_n$ are nonsingular.

The equation above can be rewritten as
$${\mathscr L}  P_n(t)  = t P_n(t)$$ with a matrix ${\mathscr L}$
such as the one
that has appeared in previous sections.

To place ourselves in the context of the {\bf bispectral problem}
we consider matrix valued polynomials $(P_n)_n$ satisfying not
only the equation above but also ``right hand side'' 
differential equations of the form:
\begin{equation}
P_nD=\Lambda_nP_n\quad\quad \text{for  all}\quad n\ge0
\end{equation}
with $\Lambda_n$ a matrix valued eigenvalue and $D$ a differential operator of order
$s$ with matrix coefficients given by
$$D=\sum_{i=0}^s \partial^iF_i(x),\quad\quad \partial=\frac{d}{dx},$$
which
acts on $P(x)$ by means of
$$P D = \sum_{i=0}^s \partial^i(P)(x)F_i(x).$$

This problem in the matrix case has been studied in
\cite{D,DG1,G10,GI,GPT1,GPT2,GPT3,GT} and in a few other places.

One can see that the differential operators that correspond to a fixed family of polynomials
form an associative algebra which in general is non-commutative, see \cite{DG1,GT}.
The problem of exhibiting elements of this algebra that have a minimal order will
occupy us in a few examples in the next two sections.

\section{A few examples}

Here we consider in detail few examples of the matrix version of the basic
Darboux process described above.

For $\lambda > 3/2$ consider the monic matrix valued polynomials which are orthogonal with respect to the weight matrix
\[
W(x) = ((2-x)x)^{\lambda-3/2} 
\begin{pmatrix} 
1 & x-1 \\ 
x-1 & 1 
\end{pmatrix},\ x \in [0,2].
\]
Let
\[
L_0 = \begin{pmatrix}
B_0 & I & & \\
A_1 & B_1 & I \\
& \ddots & \ddots & \ddots
\end{pmatrix}
\]
be the corresponding block tridiagonal matrix with
\begin{align*}
B_n &= \frac {1}{2} \frac {\lambda - 1}{(n+\lambda)(n+\lambda-1)} S + I \\
A_n &= \frac {n(n+2\lambda-2)}{4(n+\lambda-1)^2} I.
\end{align*}
Here $S = 
\begin{pmatrix}
0 & 1 \\
1 & 0 
\end{pmatrix}$ and $I = 
\begin{pmatrix}
1 & 0 \\
0 & 1 
\end{pmatrix}$.

These polynomials can be seen to be joint eigenfunctions of a first order differential operator, an observation that was made for a special value of $\lambda$ in \cite{CG1,CG2}.

If $\alpha_0$ is an arbitrary matrix we can consider the monic polynomials that result from one application of the Darboux process to the block tridiagonal matrix $L_0$ with free parameter $\alpha_0$.

We can see that for an invertible symmetric $\alpha_0$ the new orthogonality weight is given by
\begin{align*}
\widetilde{W}(x) &= (2-x)^{\lambda-3/2}x^{\lambda-5/2}
\begin{pmatrix}
1 & x-1 \\
x-1 & 1
\end{pmatrix} \\
&- \frac {-2^{2\lambda}Be\left( \frac {2\lambda-1}{2}, \frac {2\lambda-1}{2}\right)}{4(2\lambda-3)} \left( 
\begin{pmatrix}
2\lambda-2 & -1 \\
-1 & 2\lambda-2
\end{pmatrix} - (2\lambda-3)\alpha_0^{-1}\right) \delta_0(x).
\end{align*}

Here $Be$ stands for the usual Beta function.

Recall that the Darboux process played an important role in getting new bispectral situations in \cite{DG}.

We show below some examples that illustrate that for appropriate values of $\lambda$ the new polynomials are joint eigenfunctions of some higher order differential operators, i.e., we get new
bispectral situations. This appears to have little to do with $\alpha_0$ being symmetric.

\bigskip
\noindent
{\bf Example 1.}
\[
\lambda = 5/2\quad \alpha_0 = 
\begin{pmatrix}
5 & 2 \\
3 & 1
\end{pmatrix} .
\]
Here we find one differential operator $D$ satisfying
\[
P_nD = \Lambda_nP_n
\]
with
\[
D = \sum_{r=0}^4 \partial^rB_r
\]
and
\begin{align*}
B_4 &= (t-2)^2t^2 
\begin{pmatrix}
1 & -1 \\
-1 & 1
\end{pmatrix} \\
B_3 &= 4(t-2)t(3t-2) 
\begin{pmatrix}
1 & -1 \\
-1 & 1
\end{pmatrix} \\
B_2 &= \frac {24}{5} t(7t-9) 
\begin{pmatrix}
1 & -1 \\
-1 & 1
\end{pmatrix} \\
B_1 &= \frac {8}{5} 
\begin{pmatrix}
5t+6 & -8t \\
-2(5t+3) & 13t
\end{pmatrix} \\
B_0 &= \begin{pmatrix}
-\frac {8\times11}{5} & -8 \\
-\frac {32}{5} & 0
\end{pmatrix}
\end{align*}
\[
\Lambda_n = \begin{pmatrix}
\frac {(n+2)(5n^3+20n^2+3n-44)}{5} & -\frac {5n^4+30n^3+43n^2-14n+40}{5} \\
-\frac {5n^4+30n^3+43n^2+2n+32}{5} & \frac {n(5n^3+30n^2+43n+26}{5} 
\end{pmatrix}.
\]
There are no operators of lower order in the algebra.

\bigskip
\noindent
{\bf Example 2.}
\[
\lambda = 7/2,\ \alpha_0 = \begin{pmatrix}
3 & -1 \\
5 & 7
\end{pmatrix} .
\]
Here the coresponding operator is given by
\[
D = \sum_{r=0}^6 \partial^rB_r
\]
with
\begin{align*}
B_6 &= \frac {(t-2)^3t^3}{15} 
\begin{pmatrix}
1 & -1 \\
-1 & 1
\end{pmatrix} \\
B_5 &= \frac {2(t-2)^2t^2(5t-4)}{5} 
\begin{pmatrix}
1 & -1 \\
-1 & 1
\end{pmatrix} \\
B_4 &= 4(t-2)t(5t^2-8t+2) 
\begin{pmatrix}
1 & -1 \\
-1 & 1
\end{pmatrix} \\
B_3 &= 16t(5t^2-12t+6) 
\begin{pmatrix}
1 & -1 \\
-1 & 1
\end{pmatrix} \\
B_2 &= \frac {16t(131t-148)}{39} 
\begin{pmatrix}
1 & -1 \\
-1 & 1
\end{pmatrix} \\
B_1 &= \frac {32}{19} 
\begin{pmatrix}
-4(2t-7) & 3(t-6) \\
(9t-28) & -2(2t-9)
\end{pmatrix} \\
B_0 &= \frac {1}{19} 
\begin{pmatrix}
-13\times32 & -3\times64 \\
11\times32 & 0
\end{pmatrix}
\end{align*}
and we have
\[
P_nD = \Lambda_nP_n
\]
with
\begin{align*}
\Lambda_n &= \left( \frac {19n^6+285n^5+1615n^4+4275n^3+2446n^2}{285} \right) 
\begin{pmatrix}
1 & -1 \\
-1 & 1
\end{pmatrix} \\
&+ \frac {1}{285} 
\begin{pmatrix}
-12480n-6240 & 10080n-2880 \\
12960n+5280 & -10560n
\end{pmatrix}.
\end{align*}

\bigskip
\noindent
{\bf Example 3.}
\[
\lambda = 9/2\quad \alpha_0 = \begin{pmatrix}
1 & 0 \\
0 & 0
\end{pmatrix}.
\]
In this case there is one operator of order eight whose corresponding $\Lambda_n$ is given below
\[
\Lambda_n = \begin{pmatrix}
(n-3)(n+6)(n+7)(n+8)\alpha_n & -(n-2)n(n+7)(n+8)\beta_n \\
-(n+1)(n+6)\gamma_n & (n-1)n(n+1)(n+10)\delta_n 
\end{pmatrix}
\]
with
\begin{align*}
\alpha_n &= n^4+10n^3+59n^2+170n+840 \\
\beta_n &= n^4+14n^3+95n^2+322n+1080 \\
\gamma_n &= n^6+21n^5+169n^4+651n^3+1198n^2+840n-20160 \\
\delta_n &= n^4+18n^3+143n^2+558n+1512 .
\end{align*}

Once again, this corresponds to the lowest order possible differential operator
in the corresponding algebra.

\section{A few Jacobi type examples}

A different avenue for exploring the similarities as well as the differences between
the use of the Darboux process in the scalar and the matrix valued case is given by the
examples in this section.

First recall that in the scalar case it follows from results in \cite{KoKo1,Z,H2,GY1} that
the polynomials orthogonal to the weight $\mu(x)$ consisting of a Jacobi density plus two 
possible delta masses of nonnegative strengths $W,V$ at the ends of the interval, i.e.,
\[
\mu(x) = (1-x)^\alpha (1+x)^\beta
+ W  \delta_1(x) +V \delta_{-1}(x)
\]
satisfy differential equations when $\alpha$ and $\beta$ satisfy certain natural integrality conditions. The simplest example is given by the so called Koorwinder polynomials, corresponding to $\alpha=\beta=0$.
If the
weight at $1$ is the only one that is present then the order is $2 \alpha + 4$. If both
delta weights are thrown in, then the order is $2 \alpha + 2\beta +6$.
The results can be obtained by an application of the Darboux process as shown in \cite{H2,GY1}.

We consider now a small colection of situations analogous to the ones above.

The weight matrices will, as before, consist of a matrix weight density plus a pair
of deltas at the end points weighted by certain matrices $W,V$, i.e., we have
\[
\widetilde{W}(x) = (1-x)^\alpha (1+x)^\beta
\begin{pmatrix}
1 & x \\
x & 1
\end{pmatrix} \\
+ W \delta_1(x)+ V \delta_{-1}(x).
\]
For the first batch of examples we will assume that $\alpha,\beta$ are both $0$.
If $V$ and $W$ coincide with the matrix $\begin{pmatrix}
a & 0 \\
0 & 0
\end{pmatrix}$ then we find two linearly independent operators of order $5$ and one of order $6$ as well as other operators of higher order. There are no other operators of
lower order.

If $V$ is the matrix $\begin{pmatrix}
a & 0 \\
0 & 0
\end{pmatrix}$ and $W$ is the matrix $\begin{pmatrix}
0 & 0 \\
0 & b
\end{pmatrix}$ or the matrix $\begin{pmatrix}
b & 0 \\
0 & 0
\end{pmatrix}$ then we find two linearly independent operators of order $6$ as well as other operators of higher order. There are no other operators of
lower order.

More generally if $V$ is the matrix $\begin{pmatrix}
a^2 & a b \\
a b  & b^2
\end{pmatrix}$ and $W$ is the matrix $\begin{pmatrix}
c^2 &  c d\\
c d  & d^2
\end{pmatrix}$ then we have the same situation as in the last example.

In general if $V$ and $W$ are of the form $\begin{pmatrix}
a & b \\
b & c
\end{pmatrix}$ and $\begin{pmatrix}
d & e \\
e & f
\end{pmatrix}$ then the lowest order operator in the algebra is just one operator of order $8$.

We come now to a different sort of examples.

Assume that $\alpha$ and $\beta$ ($>-1$) are arbitrary, but insist in picking $W$ and
$V$ to be arbitrary (and non-necessarily equal) nonnegative multiples of the matrix $\begin{pmatrix}
1 & 1 \\
1 & 1 
\end{pmatrix}$.

In this case there is a very nice second order differential operator in the algebra which is independent of the choice of the scalar factors that appear in front
of the matrix above to give $W$ and $V$. There is no lower order operator in the
algebra.
When the deltas are both missing then
the algebra contains an operator of order $1$.

The right handed differential operator alluded to above is a scalar operator of the
usual Jacobi type, with coefficients $(1-x^2)$ and $(\alpha+\beta-1)-x(\alpha+\beta+3)$ multiplied on the right by the matrix $\begin{pmatrix}
1 & -1 \\
-1 & 1 
\end{pmatrix}$.

The eigenvalue is $-n(n+\alpha+\beta+2)$ multiplied by this same matrix.

\section{An explicit differential operator}

The paper \cite{DG} contains a proof that in the continuous-continuous case, when all operators in
question are differential operators, the so called $\ad$-conditions $$\ad L^{m+1} ( \Theta)= 0$$ are
necessary and sufficient to have what has eventually been called a bispectral situation, i.e., a solution of the
original problem tackled with Hans.

This condition gives a set of nonlinear equations that need to be solved in the unknowns $L, \Theta$.

It is important to see that this condition can be easily adapted to other situations, including the
present noncommutative one. This approach was taken up in \cite{GI}  and in \cite{GT}. In the second of these
papers the ``$\ad$-conditions'' are shown to be equivalent, once again, to bispectrality.

In general finding the differential operators of lowest possible order that appear
in a bispectral situation is not easy. By repeated applications of the Darboux
process one may end getting elements of the corresponding algebra that are not
necessarily of the lowest possible order. This issue has surfaced in several different 
papers, starting with \cite{DG} and a nice account is given in \cite{H2}.

In \cite{GT} one finds an explicit construction of a differential operator that results
form the conditions $$\ad L^{m+1} ( \Lambda) = 0.$$

The operator $D$  is given by
\[
D = \sum_{r=0}^m \partial^r(P) \frac {S_{m-r}}{r!}
\]
with matrix coefficients  $S_k=S_k(x)$ given by
\[
S_k = ((L-xI)^{m-k} \Lambda P)_0.
\]

In particular, we display some of the coefficients,  
\[
S_0 = ((L-xI)^m\Lambda P)_0,
\]
and, at the other end,
\[
S_m = (\Lambda P)_0= \Lambda_0,
\]
and this operator satisfies the desired condition 
\[
P_nD = \Lambda_nP_n \quad n \ge 0.
\]

The subindex $0$ above refers to the first entry of the corresponding ``vector'' with matrix valued entries.

\section{Toda flows with matrix valued time}

As seen in \cite{GH2}, and certainly in other places too, repeated applications of
the scalar Darboux process introduces ``times and flows'' that are related to the Toda flows.
Since these times appear as the free parameters in each application of the process, it is
only natural to raise the issue of ``matrix valued times''.

\section{Electrostatics: Heine, Stieltjes, Darboux}

In a remarkable paper that follows on earlier work of Heine, Steiltjes came 
up with a nice electrostatic interpretation for the zeros of the Jacobi
polynomials. Later work of Stieltjes as well as of I. Schur and D. Hilbert showed similar
interpretations in the case of the Laguerre and Hermite polynomials.

In \cite{G7,G8} I raise the possibility of some relation between the Darboux
process where the orthogonality functional gets more and more complicated with
every application of the process and the corresponding electrostatic interpretation of the families of polynomials that appear along the way.

It would be interesting to see what if anything of this picture
can be developed in the matrix valued case.

\section{Markov chains}

In \cite{G4,G5,G6,GdI} one finds examples of interesting quasi-birth-and-death processes that can
be studied by exploting their connection with certain specific examples of matrix valued orthogonal
polynomials. In particular in \cite{G5,G6} one finds examples where the recurrence of the process
is related to the presence of a matrix valued delta weight at $1$. Since the appearance of these
delta weights is one of the main characteristic
of an application of the Darboux process one may wonder about a probabilistic interpretation of the
relation that may exist between two Markov chains whose transition probability matrices are related
by a Darboux transformation.

\section{Things that appear before their time}

One of the most surprising phenomena uncovered in \cite{DG} has to do with
what was called ``the cusps'', namely degenerate situations that corespond to
degeneracies of ``higher order operators'' yielding ``lower order'' ones. 
To put this in the context of scalar valued orthogonal polynomials consider the
so called Koornwinder polynomials which are orthogonal to Lebesgue measure
in $[-1,1]$ plus a pair of delta masses at the end points of the interval.
In this case one knows that the corresponding orthogonal polynomials are the
common eigenfunctions of a sixth order differential operator.

In the special case when the strength of the two delta masses agrees one gets
an operator of order four, and one can say that in a search according to the order
of these operators this example, just as ``the cusps'' in \cite{DG}, appears before
its time.

We made a tentative exploration of the situation in the matrix valued case and examples
of this phenomenon are seen in section $7$.

\section{The multivariable case}

In this section we mention that in 
\cite{G9} one finds a specific random walk introduced by Hoare and Rahman, see \cite{HR}
which we show leads to a bispectral situation in terms of polynomials of two variables.
I also want to mention that in the multivariable case one finds a version of the Darboux
process to obtain interesting deformations of the two dimensional Chebyshev measure, see \cite{GI1}.

\section{Conclusion}

It is clear that very little of the development that I have tried to summarize here could have happened were it not for my good fortune of teaming up with Hans at the beginning of this journey.

\end{document}